\documentclass[12pt]{amsart}

\newtheorem{theorem}{Theorem}[section]

\newtheorem{corollary}[theorem]{Corollary}

\newtheorem{lemma}[theorem]{Lemma}

\theoremstyle{definition}

\newtheorem{remark}[theorem]{Remark}

\theoremstyle{parrafo}

\begin{document}

\title[]{The best constant for the
centered maximal operator  on  radial
decreasing functions}

\author{J. M. Aldaz and J. P\'erez L\'azaro}
\address{Departamento de Matem\'aticas,
Universidad  Aut\'onoma de Madrid, Cantoblanco 28049, Madrid, Spain.}
\email{jesus.munarriz@uam.es}
\address{Departamento de Matem\'aticas y Computaci\'on,
Universidad  de La Rioja, 26004 Logro\~no, La Rioja, Spain.}
\email{javier.perezl@unirioja.es}

\thanks{2000 {\em Mathematical Subject Classification.} 42B25}

\thanks{The authors were partially supported by Grant MTM2009-12740-C03-03 of the
D.G.I. of Spain}

%\thanks{Supported by DGICYT 1317 DP Spain}

%\subjclass{}

%\keywords{}

%\date{}

%\dedicatory{}

%\commby{}

%%% ----------------------------------------------------------------------
\begin{abstract} We show that the lowest constant appearing
in the weak type (1,1) inequality
satisfied by the centered Hardy-Littlewood maximal operator
on radial, radially decreasing integrable functions is 1.
\end{abstract}

%%% ----------------------------------------------------------------------

\maketitle

%%% ----------------------------------------------------------------------

\section {Introduction}

\markboth{J. M. Aldaz and J. P\'erez-L\'azaro}{The
maximal operator  on  radial functions}

A considerable amount of work has been devoted in the
literature to finding good bounds, or best bounds if possible,
in the inequalities satisfied by the several variants of the Hardy-Littlewood
maximal operator. We mention, for instance, \cite{A1}, \cite{A2}, \cite{A3}, \cite{A4}, \cite{A5}, \cite{ACP}, \cite{AlPe}, \cite{AlPe2}, \cite{AV}, \cite{Bou1}, \cite{Bou2}, \cite{Bou3}, \cite{Ca}, \cite{CF}, \cite{CLM}, \cite{GK}, \cite{GM}, \cite{GMM}, \cite{Ki}, \cite{Me1}, \cite{Me2}, \cite{Mu}, \cite{St1}, \cite{St2}, \cite{St3}, \cite{StSt}.
Additional references can be found in the aforementioned papers.

Let $M_d$ be the centered maximal operator (cf. (\ref{HLMF}) below for the
definition) associated to euclidean
balls and Lebesgue measure on $\mathbb R^d$. It is well known that if $1 < p \le \infty$, then there exists
a  constant $c_{p, d}$ such that for all $f\in L^p (\mathbb R^d)$, $\|M_df\|_p\le c_{p,d} \|f\|_p$, and the problem lies in determining the lowest
such
$c_{p, d}$.
When $p =\infty$, trivially we can take $c_{p, d} = 1$ for all $d$, since
averages never exceed a supremum, while if
$p=1$, then $M_d f\notin L^1(\mathbb R^d)$ unless $f = 0$ almost everywhere.
So for $p=1$ one considers instead the best constant $c_{ d}$ appearing in the
weak type $(1,1)$ inequality (cf. \ref{weaktype} below). Obviously, $c_{ d} \ge 1$, since by the Lebesgue Differentiation
Theorem, $M_d f \ge |f|$ a.e. whenever $f\in L^1 (\mathbb R^d)$.
We shall see that if we impose on $f$ the additional conditions of being
radial and radially decreasing, then actually $c_{ d} = 1$ for every dimension $d$. This improves on the previously known upper bound $c_{ d} \le 4$ (which nevertheless holds for all radial functions, not
necessarily decreasing, cf.
\cite[Theorem 3]{MS}).

Our result is obtained by
identifying the extremal case: For the class of radial, radially decreasing functions $f$
of norm one, the Dirac delta ``function" $\delta$ is extremal. That is,
$M_d \delta (x) \ge M_d f(x)$ for every $x$. Since $M_d \delta$ can be easily computed, and it yields a best constant equal to 1, the result
follows. Without the decreasing assumption on $f$, the value of the
best constant is not known (as indicated above, it is bounded
by 4).

Regarding the dependency of
$c_{p, d}$ on $d$, for general functions,   E. M. Stein showed that
 when $1 < p < \infty$, the constants $c_{p, d} < \infty$ could be chosen to be uniform in $d$ (\cite{St1}, \cite{St2}, see also \cite{St3}).

 With respect to  weak type (1,1) bounds,
best constants grow at most like $O(d)$ (cf. \cite{StSt}).
Additionally,
for the maximal
function associated to {\it cubes} (rather than euclidean balls) it is known that best bounds approach
infinity with the dimension (cf. \cite{A5}, and also \cite{Au}, where it
is shown that bounds increase at least as $O(\log^{1 - \varepsilon} d)$, for
arbitrary $\varepsilon > 0$).  While the corresponding problem for
euclidean balls has not yet been solved, it seems very likely that
uniform bounds do not exist in this case either. Hence the renewed interest
in finding natural subspaces of $L^1 (\mathbb R^d)$
for which bounds independent of $d$ can be obtained, and when uniform
bounds are known, in determining their optimal values.

\section
{Notation and results}

Let $\lambda^d$ denote the Lebesgue measure on $\mathbb R^d$, and
let $B(x, r)$ be the euclidean {\em closed} ball centered at $x$ of radius
$r>0$. Thus, $B(x, r)$ is defined using the $\ell_2$ distance
$\|x\|_2 := \sqrt{x^2_1 + \dots +x_d^2}$. The centered
maximal function $M_{d} f$ of a locally integrable
function $f$ is
\begin{equation}\label{HLMF}
M_{d} f(x) := \sup _{r > 0} \frac{1}{\lambda^d
(B(x, r))} \int _{B(x, r)} \vert f\vert d\lambda^d
\end{equation}
(the choice of closed balls in the
definition is mere convenience; using open balls instead does not
change the value of $ M_{d} f (x)$).
Likewise, the centered
maximal function $M_{d} \mu$ of a locally finite measure
$\mu$ is
\begin{equation}\label{HLMFformeas}
M_{d} \mu(x) := \sup _{r > 0} \frac{\mu (B(x, r))}{\lambda^d
(B(x, r))}.
\end{equation}

It is well known that the maximal function satisfies the following
weak type $(1,1)$ inequality:
\begin{equation}\label{weaktype}
\lambda^d(\{M_d f \ge \alpha\}) \le \frac{c_d \|f\|_1}{\alpha},
\end{equation}
where $c_d$ does not depend on $f\in L^1 (\mathbb R^d)$.

We denote the average of the
function $h$ over the
set $E$   by
\begin{equation}\label{average}
h_E:=  \frac{1}{\lambda^d
(E)} \int _{E} h d\lambda^d.
\end{equation}
Likewise, the average (with respect to Lebesgue measure)
of the
measure $\mu$ over the
set $E$  is denoted by
\begin{equation}\label{averagem}
\mu_E:=  \frac{\mu(E)}{\lambda^d
(E)}.
\end{equation}

The next ``geometric lemma on averages", states the  intuitively plausible   fact that for a radial
decreasing function on $\mathbb R^d$, the average over any ball $B$ centered at zero is
at least as large as the average over any other ball with center outside
$B$ (or on its border). By decreasing we mean non-strictly decreasing.

\begin{lemma}\label{lemma} Let $f:(0,\infty)\to (0,\infty)$ be a
decreasing function. Define
$g:\mathbb R^d\to \mathbb R$ by setting $g(x) := f(\|x\|_2)$. If
$g$ is locally integrable, then
for every pair of radii $R, r>0$, and every $y\in \mathbb R^d$ with
$\|y\|_2 \ge R$, we have
\begin{equation}\label{comparison}
g_{B(0,R)}\ge g_{B(y, r)}.
\end{equation}
\end{lemma}

\begin{remark} Actually, for the application below we only
need the case $r < R$, but since the result is also true
when $R\le r$, we do not exclude this from the statement of the lemma.
\end{remark}

\begin{remark} Obviously, if $f$ is locally integrable then so is $g$.   Local integrability of $g$ is all we need, so we only assume this weaker
condition.
\end{remark}

\begin{remark} It is natural to ask whether the hypothesis that $y$
does not belong to the interior of $B(0, R)$ can be relaxed to
$B(y, r)\setminus B(0,R) \neq \emptyset$. In fact, it is easy to
see that the latter condition is not enough, even in one dimension: Let
$\psi (x) := (1 - |x|)_+$ be the positive part of $1 - |x|$. Then $\psi_{[-1,1]} < \psi_{[-1/2,1]}$, so if
$\varepsilon > 0$ is sufficiently small, we also have $\psi_{[-1,1]} < \psi_{[-1/2,1 + \varepsilon]}$.
\end{remark}

\begin{remark} In order to
obtain large averages, one must integrate over the parts of the
space where a function is large. And this is so no matter what measure
is used.
Thus, it is tempting to conjecture that
Lemma \ref{lemma} actually holds for a large class of measures, rather
than just Lebesgue's. While this may be the case, some condition
on the measure is needed, as the following example shows.

Let $d = 2$
and set $\mu (A) := \lambda^2 (A\cap B(0,1))$, i.e., $\mu$ is
 the restriction of planar Lebesgue measure to the unit ball. Let
$\psi (x) := (1 - \|x\|_2)_+$, and observe that
$\psi_{B(0,1)} < \psi_{B(e_1, 1)}$.  This is so since
$\psi_{B(0,1)}$ is exactly equal to the average over the cone $C$ contained in $B(e_1, 1) \cap B(0,1)$ and between the lines $y = \pm \sqrt{3} x$, while
obviously  $\psi_{C} < \psi_{B(e_1, 1) \cap B(0,1)}$. But
$\mu (B(e_1, 1) \cap B(0,1)) = \mu (B(e_1, 1))$, so
$\psi_{B(e_1, 1) \cap B(0,1)} = \psi_{B(e_1, 1)}$. Therefore, Lemma \ref{lemma}
does not extend to all radial measures.
\end{remark}

\begin{remark} Another natural  attempt to generalize Lemma \ref{lemma}
 is to consider norms different from the euclidean one. Since the most
often used maximal functions on $\mathbb R^d$ are defined either using euclidean balls
or cubes with sides parallel to the axes, i.e., $\ell_\infty$ balls, the case of the
$\ell_\infty$ norm $\|x\|_\infty:= \max\{|x_1|,\dots, |x_d|\}$ is particularly interesting.
Here a function $g:\mathbb R^d\to \mathbb R$ is radial decreasing if there exists a
decreasing function $f:(0,\infty)\to (0,\infty)$ such that $g(x) = f(\|x\|_\infty)$.
Likewise, the corresponding maximal operator is obtained by averaging over $\ell_\infty$ balls $B_\infty$
in (\ref{HLMF}), instead of using $\ell_2$ balls.

The following example, due to Professor Guillaume Aubrun and included here with his permission,  shows that Lemma \ref{lemma} fails when $\|\cdot\|_2$ is replaced by $\|\cdot\|_\infty$. Let $g=\chi_{[- 2^{-1}, 2^{-1}]^d}$,
let $B_\infty(0, R) = [- 4^{-1} 3, 4^{-1} 3]^d$, and let
$B_\infty(y, r) = [4^{-1}, 4^{-1} 5] \times [- 2^{-1}, 2^{-1}]^{d-1}$,
i.e., $y = 3 e_1/4$ and $r = 1/2$. Since $2^{d + 2 } < 3^d$ provided
$d\ge 4$, for
every $d = 4, 5,\dots$ we have
$$
g_{B_\infty (0,R)}= \frac{2^d}{3^d} < \frac{1}{4} = g_{B_\infty (y, r)}.
$$
The preceding inequality is strict, so sufficiently small perturbations
of the sets involved will preserve it. Thus, Lemma \ref{lemma} also fails
 when $\|\cdot\|_2$ is replaced by $\|\cdot\|_p$, provided that
$p$ is high enough (perhaps depending on $d$). Here
$\|x\|_p := \left(|x_1|^p + \dots +|x_d|^p\right)^{1/p}$.
\end{remark}

{\em Proof of Lemma \ref{lemma}.} Observe first that the result  for $\|y\|_2 \ge R$ can be
immediately derived from the special case $\|y\|_2 = R$. To
see why, assume it holds for $\|y\|_2 = R$, and suppose
$\|w\|_2 > R$. Then
$$
g_{B(w,r)}\le g_{B(0, \|w\|_2)} \le g_{B(0,R)},
$$
since the average over a ball centered at $0$ of a radial decreasing
function does not decrease when we reduce the radius. So we assume that
$\|y\|_2 = R$. Using a change of variables if necessary, we suppose that
$R=1$ (just to simplify expressions). Then we take
$y= e_1$, by symmetry; finally, we suppose that
$f$ is left continuous. This last assumption is made purely for notational convenience: It  entails that
nonempty level sets $\{g \ge m\}$ are closed balls, agreeing with
our notation $B(0,t)$.

We show that
$r^d\int_{B(0,1)} g \ge \int_{B(e_1,r)} g$. To this end,
 it is enough to prove that for every $m>0$ the corresponding
level sets satisfy
\begin{equation}\label{levelset}
r^d\lambda^d(B(0,1)\cap \{g \ge m\}) \ge
\lambda^d(B(e_1,r)\cap \{g \ge m\}).
\end{equation}
If either $m<g(e_1)$, or $m>g(x)$ for all $x\neq 0$, then inequality (\ref{levelset}) holds trivially.
If $g(e_1) \le m \le g(x)$ for some $x\neq 0$, then  there exists a
$t\in (0,1]$ such that $\{g\ge m\}=B(0,t)$,  so it suffices to show that
\begin{equation}\label{bolas}
  r^d\lambda^d(B(0,t)) \ge
\lambda^d(B(e_1,r)\cap B(0,t)).
\end{equation}
 We assume that $r<1$ (for otherwise (\ref{bolas}) is obvious) and also that $t+r>1$ (for otherwise $B(e_1,r)\cap B(0,t)$ is the either the empty set or just one point).
With these assumptions, the boundaries of the balls $B(e_1,r)$ and $B(0,t)$ are $d-1$ spheres
whose intersection is a $d - 2$ sphere $S$, with center $c e_1$
for some $c\in (0,1)$, and
radius $\rho$. Since $B(e_1,r)\cap B(0,t)\subset B(ce_1,\rho)$, all we need
to do is to prove that $\rho \le rt$, from which (\ref{bolas}) follows.

 Let us write
$x=(x_1,\dots,x_d)$. Using symmetry, the center and the radius of the sphere $S$ can be determined by considering the intersection of $S$ with the $x_1x_2$-plane,
that is, by simultaneously solving $x_1^2 + x_2^2 = t^2$ and
$(x_1 - 1)^2 + x_2^2 = r^2$. Solving for $x_1$ yields $c =(1+t^2-r^2)/2$,
and solving for $x_2^2$, together with some elementary algebraic manipulations, gives
\begin{equation}
   \label{radius}
\rho^2 = t^2-\frac{(1 + t^2-r^2)^2}{4} = (rt)^2-\frac{(1-t^2-r^2)^2}{4} \le (rt)^2.
 \end{equation}
\qed

\begin{theorem}\label{theorem} Let $g\in L^1(\mathbb{R}^d)$ be a
 radial decreasing function. Then for every $\alpha > 0$,
\begin{equation}\label{weak}
\lambda^d(\{M_d g > \alpha\})\le\frac{\|g\|_1}{\alpha}.
\end{equation}
\end{theorem}

 \begin{proof}
Suppose $\|g\|_1 \ne 0$; using
the $1$-homogeneity of the maximal operator $M_d$ we see
that $\{M_d g > \alpha\} = \{M_d (g/\|g\|_1) > \alpha/\|g\|_1\}$, so we can always
 replace $g$ with
$g/\|g\|_1$. Thus, we assume from the start that $\|g\|_1 = 1$.
Let  $\delta$ denote the Dirac delta mass placed
 at the origin, i.e.,  $\delta$ is
the probability measure defined by $\delta(\{0\}) = 1$. In this
case it is easy to compute $M_d\delta$ explicitly: $M_d\delta (x) =
1/\lambda^d(B(x, \|x\|_2))$. Hence, for every $\alpha > 0$ the set
$\{M_d \delta \ge \alpha\}$ is a ball, and
\begin{equation}\label{weakonedelta}
\lambda^d(\{M_d \delta \ge \alpha\}) = \frac{1}{\alpha}.
\end{equation}
Inequality (\ref{weak}) is implied by (\ref{weakonedelta}), since $\delta$ is extremal in the following sense:
For every $x\in \mathbb{R}^d$ we have $M_d g(x) \le M_d \delta (x)$.
To see why, note that if $\|x\|_2 = R > 0$ and $r > 0$ is any radius, by Lemma \ref{lemma} we have
 $$
 g_{B(x, r)}\le g_{B(0, R)} \le 1/ \lambda^d(B(0, R)) = 1/\lambda^d( B(x, R)) = \delta_{B(x, R)}
 $$
 (of course, if $r\ge R$ we do not need the Lemma, since then $0\in B(x, r)$ and therefore $g_{B(x, r)}\le \delta_{B(x, r)} \le
  \delta_{B(x, R)}$). By taking the supremum over $r > 0$ we conclude that
 $M_d g(x) \le M_d \delta (x)$, as was to be shown.
 \end{proof}

 Using the preceding bound we obtain refined estimates for the operator
 norm of $M_d$, from the space of radial decreasing functions in $L^p(\mathbb{R}^d)$ to $L^p(\mathbb{R}^d)$, for
 $1 < p < \infty$. The proof, a standard Marcinkiewicz interpolation type argument, is omitted (cf., for instance \cite[p. 14]{St3}).

 \begin{corollary}\label{coro} Let $ p >1$ and let $g\in L^p(\mathbb{R}^d)$ be a radial
decreasing  function. Then
\begin{equation}\label{strong}
\|M_d g\|_p \le 2 \left(\frac{p}{p -1}\right)^{1/p} \|g\|_p.
\end{equation}
\end{corollary}

\begin{remark} Denote by $c_{p,d}$ the operator
 norm of $M_d$, from the radial decreasing functions in $L^p(\mathbb{R}^d)$ to $L^p(\mathbb{R}^d)$. Here $d$ is fixed. The preceding result entails that
 $c_{p,d}= O\left(\frac{p}{p-1}\right)$ as $p\downarrow 1$. Next we
 show by example that actually $c_{p,d}= \Theta\left(\frac{p}{p-1}\right)$,
 where $\Theta$ denotes the exact order as $p\downarrow 1$. Let $f$ be the characteristic
 function of the unit ball. Then $\|f\|_p = (\lambda^d(B(0, 1)))^{1/p}$.
 On $B(0, 1)$ the maximal function is identically one, while off
 $B(0, 1)$, writing $r = \|x\|_2$, we have $M_d f(x) \ge (r + 1)^{-d}
 \ge (2 r)^{-d}$, where the first inequality is obtained
 by averaging over the smallest ball
 centered at $x$ that fully contains $B(0,1)$, and the second
 inequality is used to
trivialize  integration in polar coordinates.
 Thus,
 \begin{equation}
 \int (M_d f)^p \ge |B(0,1)|\left(1 + \frac{1}{(p-1) 2^{d p}}\right),
 \end{equation}
 so
 \begin{equation}
c_{p,d} \ge  \left(1 + \frac{1}{(p-1) 2^{d p}}\right)^{1/p}.
 \end{equation}
Since
 $\left(1 + \frac{1}{(p-1) 2^{d p}}\right)^{1/p}\left(\frac{p - 1}{p}\right) \to 2^{-d}$
 as $p\downarrow 1$, the assertion about the exact order of $c_{p,d}$ follows.
\end{remark}

\begin{remark} If Lebesgue measure in dimension $d$ is replaced by the
standard gaussian measure, or more generally, by any finite measure
defined by a bounded, radial, radially decreasing
density, the situation is very different: The same example (one delta placed at the origin)
shows that constants for the weak type $(1,1)$ inequality
grow exponentially fast with the dimension (cf. \cite{A4})
rather than being uniformly bounded by 1. In fact,
exponential growth can be shown to hold
for some (sufficiently small) values
of $p > 1$, simply by using, instead of $\delta_0$, the characteristic function of a small
ball centered at $0$, and then arguing as in
\cite{A4}; a step in this direction is carried, for weak type
$(p,p)$ inequalities, in \cite{AlPe3}.
\end{remark}

\end{document}